\documentclass[10pt,a4paper]{article}

\usepackage{setspace}
\usepackage[utf8x]{inputenc}
\usepackage[T1]{fontenc}

\usepackage{helvet}

\usepackage{amsthm}
\usepackage{amsmath}
\usepackage{verbatim}
\usepackage{amssymb} 
\usepackage{resizegather}
\usepackage{tikz-cd}
\usepackage{babel} 
\usepackage{mathrsfs}
\usepackage{graphicx}
\usepackage{epigraph}

\usepackage{pgfplots}

\usepgfplotslibrary{fillbetween}
\pgfplotsset{compat=1.13}

\usepackage[bottom]{footmisc}

\usepackage{hyperref}

\newtheorem{theorem}{Theorem}
\newtheorem{proposition}[theorem]{Proposition}

\newtheorem{remark}[theorem]{Remark}

\newtheorem{definition}[theorem]{Definition}

\title{Motivic Milnor Fibres in Families of Real Singularities}
\author{Lars Andersen}

\begin{document}

\maketitle
\begin{abstract} In this article we prove two results concerning the motivic Milnor fibres $S^{\epsilon}(f)$ of symbol $\epsilon\in\{\pm 1, <, >\}$ associated to a map germ $f: (\mathbb{R}^n,0)\to(\mathbb{R},0)$, as defined by G. Comte and G. Fichou in \cite{CF}. Firstly, we prove that if $f,g:(\mathbb{R}^n,0)\to(\mathbb{R},0)$ are arc-analytically equivalent germs of Nash functions (\cite[Definition 7.1]{Ca}) then the virtual Poincaré polynomial of the motivic Milnor fibres $S^{\epsilon}(f)$ and $S^{\epsilon}(g)$ are equal. This extends and provides a new proof of the result \cite[Theorem 4.9]{Fi}. Secondly, let $T\subset\mathbb{R}^m$ be a real algebraic set and $f: T\times(\mathbb{R}^n,0)\to\mathbb{R}$ a polynomial function of polynomial map germs such that $f(t,0)=0$ for any $t\in T$. Then we prove that there exists a locally finite real analytic stratification $\mathcal{S}$ of $T$ such that if $S\in\mathcal{S}$ is a stratum then $f_t\sim_a f_{t'}$ are arc-analytically equivalent, for any $t,t'\in S$. Furthermore, if $T$ is compact then the stratification $\mathcal{S}$ can be taken to be finite. These two results imply in particular that the virtual Poincaré polynomials of the respective zeta functions are equal i.e $\beta(S^{\epsilon}(f_t))=\beta(S^{\epsilon}(f_{t'}))$. 
\end{abstract}
\begin{center}
    \text{Classification: }\textbf{14-XX}
\end{center}


\section{Introduction}
\subsection{Introduction}
The association of a \emph{motivic} Milnor fibre to a germ of $\mathbb{K}$-analytic map $(\mathbb{K}^n, 0)\to (\mathbb{K}, 0)$ for $\mathbb{K}$ either the field of complex or real numbers has been a subject of research for over 20 years. At the heart of the matter lies the Denef-Loeser formula \cite{DL} and its variants, which hint at a deeper connection between on the one hand the \emph{arc space} of a singularity, defined scheme-theoretically, and on the other hand its \emph{desingularisation}, a measure of the algorithmic complexity of the singularity.\footnote{this is called the \emph{Nash problem} after the pioneer work \cite{Nasharc} of John Nash dating from 1968 but left unpublished for a long time.}\\

There is a plethora of different motivic Milnor fibres. The ones we will handle concerns real analytic maps and were defined by G. Comte and G. Fichou in \cite{CF}. They are defined by way of \emph{basic semialgebraic formulas} and satisfy a Denef-Loeser formula, as proven in \cite{CF}. Yet, in contrast to the case $\mathbb{K}=\mathbb{C}$ above, no satisfactory theory of \emph{motivic integration} is yet developed for real analytic maps. And as a consequence this subject will not be discussed here since we concern ourselves here foremost with real singularities.\\

This paper has as its purpose an analysis of the motivic Milnor fibres of a semialgebraic family of polynomial functions parametrised by a real algebraic set. This will require a discussion of the arc-analytic equivalence, first defined by J.-B. Campesato in \cite{Ca}. The proof of the main theorem on the other hand, is largely based on the work \cite{PP} of A. Parusiński and P. Paunescu. Their work will be however only discussed in brief and we refer the reader to the original article for details. The notions of $\mathcal{AS}$-sets and of arc-analytic maps will likewise not be treated here and we refer to article \cite{Kurdyka1988} as well as the exposé \cite{kurdyka}.

\subsection{Outline of the Article}
We have divided the article in two parts. The first part (sections \hyperref[second part]{\ref*{second part}}-\hyperref[third part]{\ref*{third part}}) are written in the style of résumé of the subject matter whereas the second part sections \hyperref[fourth part]{\ref*{fourth part}}-\hyperref[fifth part]{\ref*{fifth part}}) consists of original research.
 
\subsection*{Acknowledgments} The author wish to express his sincere gratitude to the\\
Laboratoire de Mathématiques at USMB and especially to his thesis supervisors Georges Comte and Michel Raibaut. The material formed in this article and its successor formed a part of the authors thesis. On behalf of all authors, the corresponding author states that there is no conflict of interest.
\subsection*{Data Sharing} Data sharing not applicable to this article as no datasets were generated or analysed during the current study.

\section{Real Motivic Milnor Fibers}\label{second part}

In this section we shall discuss motivic zeta functions in the case of real isolated singularities. To put this into perspective we first recall a formula of Norbert A'Campo and briefly discuss the classical motivic zeta functions of J. Denef and F. Loeser.

\subsection{A'Campo's Formula}
The classical formula \cite[Théorème 1]{A'C} of A'Campo implies that in the case of hypersurface singularities the Lefschetz numbers of the iterates of the monodromy acting on the complex Milnor fibre can be computed from a resolution of singularities and in particular this yields a formula for the Euler characteristics of the Milnor fibre.\\

The situation is as follows. Let $f: \mathbb{C}^{n+1}\to \mathbb{C}$ be a polynomial map with $f(0)=0$ and suppose that the origin is a critical point. Write $\mathcal{F}_{\eta}=f^{-1}(\eta)\cap\bar{\mathbb{B}}_{\delta}$ for a representative of the Milnor fibre at the origin. By the resolution of singularities theorem of Hironaka \cite{hironaka} there exists a proper birational map $\sigma: M \to \mathbb{C}^{n+1}$ from a smooth connected complex algebraic variety $M$ such that the following hold.\\ 
\begin{enumerate}
\item The restriction $\sigma: M\setminus (f\circ \sigma)^{-1}(0)\to\mathbb{C}^{n+1}\setminus f^{-1}(0)$ is an algebraic isomorphism.
\item The reduced divisor $Y:=(f\circ\sigma)_{red}^{-1}(0)\subset M$ has normal crossings.
\item The reduced divisor $E:=(\sigma^{-1}(0))_{red}$ is a union of irreducible components of $Y$.
\end{enumerate}

It follows that if $Y=\bigcup_{j\in\mathcal{J}} E_j$ denotes the decomposition into irreducible components of $Y$ then there exists a subset of indices $\mathcal{K}\subset\mathcal{J}$ such that $E=\bigcup_{i\in\mathcal{K}} E_i$ is a union of irreducible components.

\begin{definition} For any $i\in\mathcal{J}$, let $N_i=\text{mult}_{E_i}(f\circ\sigma)$ be the multiplicity of $(f\circ\sigma)^{-1}(0)$ along $E_i$. For any $I\subset\mathcal{J}$ let $N_I=\gcd_{i\in I} N_i$.
\end{definition}

Let us stratify $Y$ in such a way that $E$ is a union of strata.\\ 

\begin{definition} For any nonempty subset $\emptyset\neq I\subset \mathcal{J}$ let
$$E_I^o=\bigcap_{i\in I} E_i\setminus \bigcup_{j\in\mathcal{J}\setminus I} E_j.$$
\end{definition}

The theorem of A'Campo now reads\\

\begin{theorem}[{\cite[Théorème 1]{A'C}}] For any $k\in\mathbb{N}$, the Lefschetz number $\Lambda(h_k)$ of the $k$-th iterate of the monodromy action $h:\mathcal{F}_{\eta}\to\mathcal{F}_{\eta}$ of the Milnor fibre satisfies
$$\Lambda(h^k)=\sum_{I\subset\mathcal \mathcal{K}, N_I| k} N_I \chi(E_I^0).$$
In particular the Euler characteristic of the Milnor fibre satisfies
$$\chi(\mathcal{F}_{\eta})=\Lambda(h^0)=\sum_{I\subset\mathcal \mathcal{K}} N_I \chi(E_I^0).$$
\end{theorem}

A'Campo proved this theorem in 1973, by using the spectral sequence of the sheaf of vanishing cycles as to reduce the calculation to the case of a normal crossings singularity, for which the Lefschetz numbers had been computed by him in a previous article.\\

Whilst constructing their theory of motivic integration, J. Denef and F. Loeser associated in their article \cite{DL} to any nonconstant morphism 
$$f: X\to\mathbb{A}^1_k$$
from an irreducible nonsingular algebraic variety over a field $k$ of characteristic zero, motivic zeta functions. These are formal power series with coefficients in a ring of motives, which is to say in this setting, a localisation of the Grothendieck ring of varieties $K_0(\text{Var}_k)$. The coefficients are then the motives corresponding to certain elements of the arc spaces introduced by Nash. A certain formal limit of the motivic zeta functions is then called the \emph{motivic Milnor fibre} of $f$. We will not delve further into this matter, lest to say that when $k=\mathbb{C}$ one recovers most of the usual invariants of the complex Milnor fibre (such as the Euler characterstic) from the motivic Milnor fibre, by applying a certain realisation morphism.\\

We will however discuss a version for real polynomial maps, due to G.Comte and G.Fichou, from which one obtains the Euler characteristic of the real Milnor fibres. As their construction is based on ``motives'' as given by a localisation of the ring of \emph{basic semialgebraic formulas}, we first need to discuss the latter.

\subsection{Basic Semialgebraic Formulas}
In the article \cite{CF} G. Comte and G. Fichou introduced the Grothendieck ring $K_0(\text{BSA}_{\mathbb{R}})$ of basic semialgebraic formulas. They were brought to do so because of the lack of algebraic structure of the Grothendieck ring $K_0(\text{SA}_{\mathbb{R}})$ of real semialgebraic sets; indeed the latter ring is isomorphic to the ring of integers generated by the class of a point. The idea then was to consider not semialgebraic sets themselves but rather the formulas defining them and thus produce a ring with richer structure. In what follows we will denote $\text{Var}_{\mathbb{R}}$ the category of varieties over $\text{Spec}(\mathbb{R})$ and if $X\in\text{Var}_{\mathbb{R}}$ is a real variety we shall denote by $X(\mathbb{R})$ its set of real points.

\begin{definition}[{\cite[§ 1.3]{CF}}] A basic semialgebraic formula in $n$ variables $A$ is a finite number of equalities, inequalities and inequations of polynomials in $n$ variables. That is, there exists $I,J, L\subset\mathbb{N}$ and polynomials $p_i,q_j,r_l\in\mathbb{R}[x_1,\dots,x_n]$ for $i\in I, j\in J$ and $l\in L$ such that
$$A=\{p_i=0,\quad q_j\neq 0,\quad r_l>0,\quad \forall i\in I, \forall j\in J, \forall l\in L\}.$$
A basic semialgebraic formula without inequalities and inequations is called algebraic. The set of semialgebraic formulas in $n$ variables is denoted by $\text{BSA}_{n}$. The union of the sets $\text{BSA}_n$ over $\mathbb{N}$ is denoted by $\text{BSA}$. 
\end{definition}

\begin{remark} A basic semialgebraic formula is not determined by its set of real points
$$A(\mathbb{R})=\{\bar{x}\in \mathbb{R}^n\ |\ p_i(\bar{x})=0,\quad q_j(\bar{x})\neq 0,\quad r_l(\bar{x})>0,\quad i\in I, j\in J, l\in L \}.$$
For instance $A=\{x^2+a>0\}$ and $B=\{x^2+b>0\}$ with $a,b$ positive real numbers, define different formulas although their sets of real points are the same.
\end{remark}

An isomorphism of formulas is defined as follows\\
\begin{definition}[{\cite[Definition 2.6]{CF}}] Let 
$$A=\{A',p_i>0, i=1,\dots,l\},\qquad B=\{B', q_i>0, i=1,\dots,l\}$$
be formulas such that $A',B'\subset\mathbb{R}^n$ are Zariski constructible subsets and $p_i, q_i\in\mathbb{R}[x_1,\dots, x_n]$. We say that $A$ and $B$ are isomorphic if there exists an isomorphism of real algebraic varieties
$$\phi: \tilde{A}\to \tilde{B}$$
where
$$\tilde{A}=\{A', p_1=z_1^2,\dots, p_l=z_l^2\},$$
$$\tilde{B}=\{B', q_1=z_1^2,\dots, q_l=z_l^2\}$$ 
such that $\phi$ is equivariant with respect to the action of $\{\pm 1\}^l$ on $\tilde{A}$ and $\tilde{B}$ and such that the morphism induced by extension of scalars $\phi_{\mathbb{C}}: \tilde{A}_{\mathbb{C}}\to\tilde{B}_{\mathbb{C}}$ is a complex algebraic isomorphism.
\end{definition}

\begin{definition}[{\cite[§ 1.3]{CF}}]
Let $G$ be the free abelian group $\langle [A]| A\in \text{BSA}\rangle$ of equivalence classes of basic semialgebraic formulas where if $A$ and $B$ are algebraic, that is have no inequalities, and if $A\cong B$ by an algebraic isomorphism then $[A]=[B]$. Let $K_0(\text{BSA}_{\mathbb{R}})$ be the quotient of $G$ by the relations
\begin{enumerate}
\item $[A, q=0]+[A, q\neq 0]=[A],\quad A\in \text{BSA}_{n,\mathbb{R}},\quad q\in\mathbb{R}[x_1,\dots,x_n],$
\item $ [A, q>0]+[A, q<0]+[A, q=0]=[A],\quad A\in\text{BSA}_{n,\mathbb{R}},q\in\mathbb{R}[x_1,\dots,x_n].$
\item If $A,B$ are basic semialgebraic formulas with disjoints sets of variables then 
$$[A][B]=[A, B]$$
where $A, B$ is the conjunction of the formulas $A$ and $B$.
\end{enumerate}

\end{definition}

Recall that given a field $k$ the Grothendieck group\footnote{More generally, in any exact category there is a notion of Grothendieck group} $K_0(\text{Var}_k)$ of algebraic varieties over the spectrum of $k$ is formed by isomorphism classes of varieties in such a way that one has the excision property for closed subvarieties. 

\begin{definition}[{\cite[Definition 1.1]{McCrory}}] Let $G$ be the free abelian group of isomorphism classes $[X]$ of varieties $X$ over the spectrum of $k$. The Grothendieck group $K_0(\text{Var}_k)$ is the quotient of $G$ modulo the relation
\begin{enumerate}
    \item If $Y\subset X$ is a closed subvariety then $[X]=[X\setminus Y]+[Y]$.
\end{enumerate}
The quotient of $G$ by the above relation and by the supplementary relation 
$$[X\times_k Y]=[X][Y]$$
gives $K_0(\text{Var}_k)$ the structure of a commutative ring. The class $[\mathbb{A}^1_{k}]$ of the affine line is denoted by $\mathbb{L}_k$ and is called the Lefschetz motive.
\end{definition}

We have that $K_0(\text{Var}_{\mathbb{R}})$ is isomorphic to a subring of $K_0(\text{BSA})$ by the following proposition:\\

\begin{proposition}[{\cite[Proposition 1.3]{CF}}] The morphism 
$$i: K_0(\text{Var}_{\mathbb{R}})\to K_0(\text{BSA}_{\mathbb{R}})$$
defined by sending the class of an algebraic variety to the formula defining it, is an injection of rings.
\end{proposition}

By the universal property of the Grothendieck ring $K_0(\text{Var}_k)$ if $(R, \cdot, +)$ is a commutative ring and if $e: \text{Var}_k\to R$ is any application such that
\begin{enumerate}
\item $e(X)=e(Y)$ if $X$ and $Y$ are isomorphic $k$-varieties,
\item $e(X)=e(Y)+e(X\setminus Y)$ if $Y\subset X$ is a closed subvariety,
\item $e(X\times_k Y)=e(X)\cdot e(Y)$
\end{enumerate}

then there exists a unique ring homomorphism $\tilde{e}: K_0(\text{Var}_k)\to R$ such that $\tilde{e}([X])=e(X)$ for any $X\in\text{Var}_k$. Applications $e: \text{Var}_k\to R$ as above are called (additive and multiplicative) invariants in the literature and ring homomorphisms $\tilde{e}: K_0(\text{Var}_{k})\to R$ are called realisations of $K_0(\text{Var}_k)$. When $k=\mathbb{C},\mathbb{R}$ the Euler characteristics with compact supports gives an example of a realisation of $K_0(\text{Var}_k)$; another is given by the virtual Poincaré polynomial\footnote{In the complex case another example is given by the Hodge-Deligne polynomial} the existence of which in the real case was first proven by C. McCrory and A. Parusi\'nski.\\

\begin{definition}[{\cite[Definition 2.3]{McCrory}}] The virtual Poincaré polynomial is the unique ring morphism $\beta_{\text{Var}}: K_0(\text{Var}_{\mathbb{R}})\to\mathbb{Z}[u]$ such that
$$\beta_{\text{Var}}(X)(u)=\sum_{i\geq 0} (-1)^{i}b_i(X)u^i$$
whenever $X$ is a compact nonsingular real variety with Betti numbers 
$$b_i(X)=\dim_{\mathbb{Z}/2\mathbb{Z}} H_i(X; \mathbb{Z}/2\mathbb{Z}),\qquad i\geq 0$$
in integral homology with $\mathbb{Z}/2\mathbb{Z}$-coefficients. 
\end{definition}

This morphism can be extended to a morphism on the Grothedieck ring of basic semialgebraic formulas $K_0(\text{BSA})$, giving a realisation of this ring.  The definition is by induction on the number of inequalities appearing in a representative of a class $[A]\in K_0(\text{BSA})$. More precisely\\

\begin{definition}[{\cite[Proposition 3.1]{CF}}]\label{def polynome} The virtual Poincaré polynomial of basic semialgebraic formulas is the unique morphism of rings $\beta: K_0(\text{BSA})\to\mathbb{Z}[u][1/2]$ such that
\begin{enumerate}
\item If $A\in \text{BSA}$ and $p\in\mathbb{R}[x_1,\dots, x_n]$ then 
 $$\beta([A, p>0])=\frac{1}{4}\beta([A, p=z^2])-\frac{1}{4}\beta([A, p=-z^2])$$
 $$+\frac{1}{2}\beta([A, p\neq 0]).$$
\item If $[A]\in K_0(\text{Var}_{\mathbb{R}})$ then $\beta([A])=\beta_{\text{Var}_{\mathbb{R}}}([A])$.
\end{enumerate}
\end{definition}

The above morphism has the following geometric property that evaluated at $-1$ one recovers the Euler characteristic with compact support of the set of real points of the formula:\\ 
\begin{proposition}[{\cite[Proposition 3.4]{CF}}] For any basic semialgebraic formula $A$,
$$\beta([A])(-1)=\chi_c(A(\mathbb{R}))$$
where $\chi_c$ denotes Euler characteristic with compact support.
\end{proposition}

One thus gets an extension\footnote{As a matter of fact, the virtual Poincaré polynomial was defined by Comte and Fichou exactly as to impose this property.} of $\chi_c: K_0(\text{Var}_{\mathbb{R}})\to\mathbb{Z}$ to a morphism 
$$\chi_c: K_0(\text{BSA})\to\mathbb{Z}[1/2].$$ 
The following proposition implies that isomorphic formulas have the same images by the virtual Poincaré polynomial.\\
\begin{proposition}\label{invariants for formulas} Suppose given formulas $A,B\in \text{BSA}_{\mathbb{R}}$. If there exists an isomorphism of formulas $\phi: A\to B$ then the classes of $A$ and $B$ have the same image $\beta([A])=\beta([B])$ under $\beta: K_0(\text{BSA}_{\mathbb{R}})\to\mathbb{Z}[u][1/2]$
\end{proposition}
\begin{proof} This follows from \cite[Proposition 2.2]{CF} and from the fact that isomorphic real algebraic varieties define the same class in $K_0(\text{Var}_{\mathbb{R}})$ hence have the same images under $\beta_{\text{Var}}$.
\end{proof}
\begin{remark}[``Topology of Formulas''] Little is known on the exact relationship between on the one hand the topology of the set $A(\mathbb{R})$ of real points of a formula and on the other hand the value $\beta([A])$ of the virtual Poincaré polynomial of the class of the formula in question. The problem is quite interesting but indeed very complicated, for, in order to find $\beta([A])$ one is required to caculate $3^k$ terms where $k$ is the number of inequalities appearing in $A$. To calculate just one term one is then required to find the Betti numbers of a real Zariski-constructible set. By additivity one is eventually reduced to the case of a compact nonsingular algebraic variety by compactifying and resolving singularities. But finding the Betti numbers, or even the zeroeth one, of a compact, nonsingular real algebraic set is a notably difficult problem, albeit that one knows by the results of Abkalut and King \cite{king} that it is diffeomorphic to some smooth manifold. 
\end{remark}

\subsection{Real Motivic Milnor Fibres}
In analogy with the motivic zeta functions of Denef and Loeser \cite{DL} for nonconstant complex algebraic morphisms, and the motivic zeta functions of S. Koike and A. Parusiński \cite{koike} and of G. Fichou \cite{Fi} for real analytic maps, Comte and Fichou associated in the article \cite{CF} motivic zeta functions to any real polynomial map $f: \mathbb{R}^{n+1}\to\mathbb{R}$ with $f(0)=0$. As already stated, these motivic zeta functions are formal power series with coefficients given by motives corresponding to certain subsets of truncated arc-spaces. By a motive is meant here an element of the localisation $K_0(\text{BSA})[\mathbb{L}_{\mathbb{R}}^{-1}]$ of the Grothendieck ring of basic semialgebraic formulas with respect to the Lefschetz motive. We recall in what follows the construction of said motivic zeta functions and their known properties.\\

Consider the \emph{arc scheme} $\mathcal{L}(\mathbb{A}_{\mathbb{R}}^{n+1})$ over the spectrum of $\mathbb{R}$. Its set of real points $\mathcal{L}(\mathbb{A}_{\mathbb{R}}^{n+1})(\mathbb{R})$ is the set of $\mathbb{R}[[t]]$-rational points of $\mathbb{A}_{\mathbb{R}}^{n+1}$. There exists a subscheme $\mathcal{L}(\mathbb{A}_{\mathbb{R}}^{n+1}, 0)$ having as set of real points
$$\mathcal{L}(\mathbb{A}_{\mathbb{R}}^{n+1},0)(\mathbb{R})=\{\gamma\in (\mathbb{R}[[t]])^{n+1}\ |\ \gamma(0)=0\}.$$
For each $k\in\mathbb{N}$ there exists a $\mathbb{R}$-scheme of finite type $\mathcal{L}_k(\mathbb{A}_{\mathbb{R}}^{n+1})$, called the \emph{scheme of truncated arcs} of order $k+1$, having as set of real points the set of $\mathbb{R}[[t]]/t^{k+1}$-rationals points of $\mathbb{A}_{\mathbb{R}}^{n+1}$.\\

In what follows we shall only consider the set of real points of these schemes, and we shall write $\mathcal{L}(\mathbb{R}^{n+1}), \mathcal{L}(\mathbb{R}^{n+1}, 0)$ and $\mathcal{L}_k(\mathbb{R}^{n+1},0)$ instead of $\mathcal{L}(\mathbb{A}_{\mathbb{R}}^{n+1})(\mathbb{R}),$\\
$\mathcal{L}(\mathbb{A}_{\mathbb{R}}^{n+1},0)(\mathbb{R})$ and $\mathcal{L}_k(\mathbb{A}_{\mathbb{R}}^{n+1},0)(\mathbb{R})$, respectively.\\

For each $k\in\mathbb{N}$ there exists a truncation map 
$$\pi_k: \mathcal{L}(\mathbb{R}^{n+1}, 0)\to \mathcal{L}_k(\mathbb{R}^{n+1}, 0)$$
sending an arc to its $(k+1)$-th truncation.  For a truncated arc  
$$\gamma_k\in\mathcal{L}_k(\mathbb{R}^{n+1},0),\qquad \gamma_k(t)=a_k t^k+a_{k-1}t^{k-1}+\dots+a_0$$ 
one defines a function
$$\text{ord}: \mathcal{L}_k(\mathbb{R}^{n+1},0)\to\mathbb{N}\cup\{\infty\}$$
sending $\gamma_k$ to $\text{ord}(\gamma_k)$ and sending $\gamma_k\equiv 0$ to $\text{ord}(0)=\infty$. There also exists a map, the angular component map
$$\text{ac}: \mathcal{L}_k(\mathbb{R}^{n+1},0)\to\mathbb{R},\qquad\text{ac}(\gamma_k)=a_{\text{ord}\gamma_k}$$
which sends a truncated arc to its leading coefficient.\\

One then defines basic semialgebraic formulas 
$$X_{k,f}^{1}:=\{\gamma\in\mathcal{L}_k(\mathbb{R}^{n+1},0)\ |\ \text{ord}(f\circ\gamma)=k,\quad \text{ac}(f\circ\gamma)=1\},$$
$$X_{k,f}^{-1}:=\{\gamma\in\mathcal{L}_k(\mathbb{R}^{n+1},0)\ |\ \text{ord}(f\circ\gamma)=k,\quad \text{ac}(f\circ\gamma)=-1\},$$
$$X_{k,f}^{>}:=\{\gamma\in\mathcal{L}_k(\mathbb{R}^{n+1},0)\ |\ \text{ord}(f\circ\gamma)=k,\quad \text{ac}(f\circ\gamma)>0\},$$
$$X_{k,f}^{<}:=\{\gamma\in\mathcal{L}_k(\mathbb{R}^{n+1},0)\ |\ \text{ord}(f\circ\gamma)=k,\quad \text{ac}(f\circ\gamma)<0\}.$$
The motivic zeta functions corresponding to $f$ are then defined as follows\\
\begin{definition}[{\cite[§ 4.1]{CF}}] For any symbol $\epsilon\in\{\pm 1, <,>\}$ define
$$Z^{\epsilon}(f):=\sum_{k\geq 0}[X_{k,f}^{\epsilon}]\mathbb{L}^{-(n+1)k}T^k\in K_0(\text{BSA}_{\mathbb{R}})[\mathbb{L}^{-1}][T].$$
\end{definition}

Note that $Z^{\pm 1}(f)$ are elements of $\mathcal{M}_{\mathbb{R}}[T]=K_0(\text{Var}_{\mathbb{R}})[\mathbb{L}^{-1}][T]$, that is, their coefficients are motives in the sense of Denef and Loeser.\\ 
\begin{remark} By applying $\beta: K_0(\text{BSA})\to \mathbb{Z}[1/2]$ to $Z^{\pm 1}(f)$ one recovers the motivic zeta functions of G. Fichou (see \cite{Fichou}). This follows from \cite[Proposition 1.3]{CF}.
\end{remark}

In order to prove the rationality of the above motivic zeta functions, Comte and Fichou established a certain formula for $Z^{\epsilon}(f)$, analoguous to the formula \cite[Theorem 2.4]{DL} established by Denef and Loeser for their motivic zeta functions. The content of this formula is that one can compute the zeta functions via a resolution of singularities. To describe it one starts with a proper birational morphism $\sigma: M\to\mathbb{R}^{n+1}$ from a smooth manifold $M$ such that\\
\begin{enumerate}
\item The restriction $\sigma: M\setminus (f\circ\sigma)^{-1}(0)\to\mathbb{C}^{n+1}\setminus f^{-1}(0)$ is an algebraic isomorphism,
\item The reduced divisors of $f\circ\sigma$ and $\det\text{Jac }\sigma$ are normal crossings,
\item The reduced divisor $E=(\sigma^{-1}(0))_{red}$ is a union of irreducible components of $Y:=((f\circ\sigma)^{-1}(0))_{red}$.\\ 
\end{enumerate}

Under these hypotheses one can write
$$Y=\bigcup_{j\in\mathcal{J}} E_j,\qquad E=\bigcup_{k\in\mathcal{K}} E_k,$$
for some subsets $\mathcal{J}\subset\mathbb{N}$ and $\mathcal{K}\subset \mathcal{J}$.\\

\begin{definition}\label{defn 10} For any $i\in \mathcal{J}$ let
$$N_i=\text{mult}_{E_i}(f\circ\sigma),\qquad \nu_i-1:=\text{mult}_{E_i}\det\text{Jac}(f\circ\sigma)$$
\end{definition}

As the divisor of $f\circ\sigma$ is normal crossings in $M$ there exists for any point $p\in Y$ an affine open $U\subset M$ containing $p$ a regular sequence of regular functions $x_i: U\to \mathbb{R}$ and a unit $u\in \mathcal{O}_U^{\times}$ such that
\begin{equation}\label{local thang}
f\circ\sigma(x)=u(x)\prod_{i\in\mathcal{J}} x_i^{N_i},\qquad\forall x\in U.
\end{equation}
One now stratifies $Y$ in such a way that $E$ is a union of strata.\\

\begin{definition}
For any nonempty subset $I\subset \mathcal{J}$ let 
$$E_{I}^o:=\bigcap_{i\in I} E_i\setminus\bigcup_{j\in \mathcal{J}\setminus I} E_j$$
and let $N_I=\gcd_{i\in I}(N_i)$.
\end{definition}

For any symbol $\epsilon\in \{\pm 1, <,>\}$ and for any stratum $E_I^o$ one constructs a basic semialgebraic formula $E_I^{o,\epsilon}$ such that $E_I^{o,1}(\mathbb{R})$ is the set of real points of an étale covering of degree $N_I$ of the extension of scalars of $E_I^o$ to $\mathbb{C}$. For this to indeed yield a basic semialgebraic formula one first has to define it locally on affine open sets, using (\hyperref[local thang]{\ref*{local thang}}).

\begin{definition}[{\cite[§ 4.1]{CF}}\label{defn 12}] For any symbol $\epsilon\in\{\pm 1, <,>\}$, if $U\subset M$ is the affine open set in (\hyperref[local thang]{\ref*{local thang}}) above then 
$$R_{I,U}^{1}:=\{(x,t)\in (E_{I}^o\cap U)\times\mathbb{R}\ |\ t^{N_I}u(x)=1\},$$
$$R_{I,U}^{-1}:=\{(x,t)\in (E_{I}^o\cap U)\times\mathbb{R}\ |\ t^{N_I}u(x)=-1\},$$
$$R_{I,U}^{>}:=\{(x,t)\in (E_{I}^o\cap U)\times\mathbb{R}\ |\ t^{N_I}u(x)>0\},$$
$$R_{I,U}^{<}:=\{(x,t)\in (E_{I}^o\cap U)\times\mathbb{R}\ |\ t^{N_I}u(x)<0\}$$
where $N_I$ is an in Definition \hyperref[defn 10]{\ref*{defn 10}} and $u$ as in (\hyperref[local thang]{\ref*{local thang}}) above. For a covering $M=(U_l)_{l\in L}$, set 
$$[E_{I}^{0,\epsilon}]=\sum_{S\subset L} (-1)^{|S|+1}[R_{I,\cap_{s\in S} U_s}^{\epsilon}]\in K_0(\text{BSA}).$$
\end{definition}

The rationality of the zeta functions $Z^{\epsilon}(f)(T)$ follows from the following theorem.\\
\begin{theorem}[{\cite[Theorem 4.2]{CF}, \cite{DL}}] For any symbol $\epsilon\in\{\pm 1, <,>\}$ the zeta functions $Z^{\epsilon}(f)(T)$ satisfy
$$Z^{\epsilon}(f)(T)=\sum_{I\cap \mathcal{K}\neq \emptyset}(\mathbb{L}-1)^{|I|-1}[\tilde{E}_{I}^{0,\epsilon}]\prod_{i\in I}\frac{\mathbb{L}^{-\nu_i}T^{N_i}}{1-\mathbb{L}^{-\nu_i}T^{N_i}}.$$
\end{theorem}

One considers the following formal limits of the zeta functions above:\\
\begin{definition}[{\cite[Definition 4.5]{CF}}]\label{motivic fibres} For any symbol $\epsilon\in\{\pm 1, <,>\}$, let 
$$S^{\epsilon}(f):=\sum_{I\cap K\neq \emptyset}(1-\mathbb{L})^{|I|-1}[\tilde{E}_{I}^{0,\epsilon}]\in K_0(\text{BSA}_{\mathbb{R}}).$$
\end{definition}

In accordance with the terminology used in the article \cite{CF} (see also \cite{DL}) one calls $S^{\pm 1}(f)$ the \emph{motivic positive and negative Milnor fibres} of $f$ at the origin. Introducing new terminology we will refer to $S^{>}(f)$ as the \emph{motivic positive Milnor tube}, and to $S^{<}(f)$ as the \emph{motivic negative Milnor tube}. They do not depend on the choice of a strong desingularisation, by \cite[Remark 4.3]{CF}.\\

In the case of isolated singularities one recovers the Euler characteristics of the classical Milnor fibres by applying the realisation  
$$\tilde{\chi}_c: K_0(\text{BSA}_{\mathbb{R}})\to\mathbb{Z}[1/2],\qquad \tilde{\chi}_c([X]):=\beta([X])(-1)$$
as according to the following real (partial) analogue of the A'Campo formula. So suppose that $f: \mathbb{R}^{n+1}\to \mathbb{R}$ is a polynomial function having an isolated critical point in the origin. Let $(\epsilon_0, \delta_0)$ be Milnor data of $f$ at the origin and for any $\delta\in (0, \delta_0)$ and any $\eta\in (0, \epsilon(\delta))$ write 
$$\mathcal{F}_{\eta}^{>}:=f^{-1}((0,\eta))\cap \mathbb{B}_{\delta},$$
$$\mathcal{F}_{\eta}^{<}:=f^{-1}((-\eta,0))\cap \mathbb{B}_{\delta},$$
and write
$$\mathcal{F}_{\eta}^{1}=f^{-1}(\eta)\cap \mathbb{B}_{\delta},$$
$$\mathcal{F}_{\eta}^{-1}=f^{-1}(-\eta)\cap \mathbb{B}_{\delta}$$
for the positive and negative Milnor fibres at the origin.

\begin{theorem}[{\cite[Theorem 4.12]{CF}}] Let 
$$f: \mathbb{R}^{n+1}\to\mathbb{R},\qquad f(0)=0$$ 
be a polynomial function having the origin as an isolated critical point. Then for any $\delta\in (0,\delta_0)$ and any $\eta\in (0,\epsilon(\delta))$,
$$\tilde{\chi}_c(S^{\epsilon}(f))=(-1)^{n}\chi_c(\mathcal{F}_{\eta}^{\epsilon}),\qquad \epsilon\in\{\pm 1,<,>\}$$
where $\chi_c$ denotes Euler characteristic with compact support.
\end{theorem}

\section{The Arc-Analytic Equivalence}\label{third part}
\subsection{Nash Functions and Manifolds}
We start by first recalling the arc-analytic equivalence for Nash germs, introduced by J-B Campesato in the article \cite{Ca}.

\begin{definition}[{\cite{nash} \cite[Definition 6.1]{Ca}}] A function $f: U\to \mathbb{R}$ defined on an open semialgebraic set $U\subset\mathbb{R}^n$ is said to be \emph{Nash} if it is semialgebraic and of class $C^{\infty}$. A map $f: U\to\mathbb{R}^m$ is said to be Nash if each of its components functions are Nash.
\end{definition}

By \cite[Proposition 8.1.8]{roy} one has that a function $f: U\subset\mathbb{R}^n\to\mathbb{R}$ is Nash if and only if it is real analytic and satisfies a nontrivial polynomial equation. In particular, by \cite[Corollary 8.1.6]{roy}, the ring of germs of Nash functions is isomorphic to the ring of algebraic power series. This yields a category whose objects are semialgebraic sets and whose arrows are Nash functions. 
The isomorphisms in this category are then: 
\begin{definition} A Nash diffeomorphism is a semialgebraic map $\phi: U\to V$ such that $\phi$ and its inverse are Nash functions.
\end{definition}

The notion of manifold in this category is given by the following.
\begin{definition}[{\cite[Definition 6.3]{Ca}}] A Nash manifold of dimension $d\in\mathbb{N}$ is a semialgebraic subset $M\subset\mathbb{R}^n$ such that for any $x\in M$ there exists an open semialgebraic neighborhood $U\subset\mathbb{R}^n$ of $x$, an open semialgebraic neighborhood $V\subset\mathbb{R}^n$ of the origin and a Nash diffeomorphism $\phi: U\to V$ with $\phi(x)=0$ such that $\phi(M\cap U)=\mathbb{R}^d\times\{0\}$.
\end{definition}

\subsection{Nash-Isomorphic $\mathcal{AS}$-Sets}
The virtual Poincaré polynomial on compact real algebraic varieties is invariant under Nash diffeomorphisms, as follows from the result \cite[Theorem 3.3]{Fi} (compare \cite[Proposition 1.2]{link}).\\

For noncompact varieties this is no longer true; one gets a counter-example thereof by taking $X\subset\mathbb{R}^2$ to be a hyperbola. Then $X$ is Nash diffeomorphic to the union $Y$ of two disjoint affine lines $\mathbb{A}^1$. So on the one hand $\beta(X)=u-1$ because if $\bar{X}\subset\mathbb{P}^2$ is a projective compactification then $\bar{X}\setminus X$ consists of two points, whereas on the other hand $\beta(Y)=2\beta(\mathbb{A}^1)=2u$ by additivity. In the case of $\mathcal{AS}$-sets one defines 

\begin{definition}[{\cite[Definition 3.2]{Fi}}] If $A,B\in\mathcal{AS}$ then $A$ and $B$ are said to be \emph{Nash-isomorphic} if there exists compact Nash manifolds $A\subset X$ and $B\subset Y$ and a Nash diffeomorphism $\phi: X\to Y$ such that $\phi(A)=B$.
\end{definition}

Then Fichou proved \cite[Theorem 3.3]{Fi} that the virtual Poincaré polynomial as defined on the $\mathcal{AS}$-collection (which contains the real algebraic sets) is invariant under Nash isomorphisms.\\

\subsection{The Arc-Analytic Equivalence}
We now turn to the arc-analytic equivalence.

\begin{definition}[{\cite[Definition 7.1]{Ca}}]\label{arc equivalence}
Let $f,g:(\mathbb{R}^n,0)\to(\mathbb{R},0)$ be germs of Nash functions. If there exists a germ of semialgebraic homeomorphism 
$$h:(\mathbb{R}^n,0)\to(\mathbb{R}^n,0)$$
with $f=g\circ h$ such that $h$ is arc-analytic and if there exists a constant $c>0$ such that
$$|\det\text{Jac}(h)|>c$$
there where the Jacobian determinant is defined, then $f$ is said to be arc-analytically equivalent to $g$ via $h$ and one writes $f\sim_a g$ or $f\sim g$ if no ambiguity is possible.
\end{definition}

This notion is due to J.-B. Campesato. He introduced it in order to prove that the \emph{Blow-Nash equivalence} on germs of Nash functions \cite{Fi}, \cite{Ca} is indeed an equivalence relation. To describe the latter we need the notion of Nash modification. To this aim, let us first recall the notion of complexification of a real analytic manifold.

\begin{definition}[{\cite[Definition 1.4]{Ca}}] A complexification of a real analytic manifold $M$ of dimension $n$ is a complex analytic manifold $M_{\mathbb{C}}$ endowed with a real analytic isomorphism 
$$\phi: M\to \phi(M)\subset M_{\mathbb{C}}$$
with $\phi(M)$ a real analytic subvariety of $M_{\mathbb{C}}$ such that for any $x\in M_{\mathbb{C}}$ there exists a neighborhood $U_{\mathbb{C}}\subset M_{\mathbb{C}}$ of $x$ and a complex analytic isomorphism 
$$\psi: U_{\mathbb{C}}\to \psi(U_{\mathbb{C}})\subset\mathbb{C}^n$$
such that
$$\psi(\phi(M)\cap U_{\mathbb{C}})\subset \mathbb{R}^n\cap\psi(U_{\mathbb{C}}).$$
\end{definition}

In other words $M$ can be identified with a real analytic subvariety of its complexification via $\phi$ and this identification respects the charts. One can show that any real analytic manifold possesses a complexification. Moreover the construction is functorial in the sense that if $g: M\to N$ is a proper real analytic morphism of manifolds and if $M_{\mathbb{C}}$ and $N_{\mathbb{C}}$ are complexifications then $g$ extends to a proper complex analytic map $g_{\mathbb{C}}: U\to V$ of open neighborhoods 
$U\subset M_{\mathbb{C}}$ and $V\subset N_{\mathbb{C}}$ of $M$ and $N$, respectively. The complexification $g_{\mathbb{C}}$ is therefore unique as a germ of complex analytic maps. 

\begin{definition}[{\cite[§ 6.2]{Ca}}] A Nash function $\sigma: U\subset\mathbb{R}^n\to\mathbb{R}^p$ is said to be a Nash modification if it is proper, surjective and if its complexification is proper and bimeromorphic.
\end{definition}

\begin{definition}[{\cite[Definition 4.2]{Fi}}]\label{defn 4.2} Let $f,g:(\mathbb{R}^n,0)\to(\mathbb{R},0)$ be germs of Nash functions. One says that $f$ and $g$ are blow-Nash equivalent if the following holds.
\begin{enumerate}
\item There exist two germs of Nash modifications 
$$\sigma_f: (M,\sigma_f^{-1}(0))\to(\mathbb{R}^n,0),\qquad \sigma_g: (M',\sigma_g^{-1}(0))\to(\mathbb{R}^n,0)$$
such that the associated divisors of $f\circ\sigma_f$ and $g\circ\sigma_g$ are normal crossings, and such that the associated divisors of $\text{Jac}(\sigma_f)$ and $\text{Jac}(\sigma_g)$ are normal crossings. 

\item There exists a germ of Nash isomorphisms 
$$\Psi: (M,\sigma_f^{-1}(0))\to(M',\sigma_g^{-1}(0))$$
preserving the multiplicities of the divisors $\text{div}\text{Jac}(\sigma_f)$ and $\text{div}\text{Jac}(\sigma_g)$ along each of the irreducible components of the exceptional divisors $\sigma_f^{-1}(0)$ and $\sigma_g^{-1}(0)$ of the modifications. 

\item The map $\Psi$ induces a germ of semialgebraic homeomorphism 
$$\psi: (\mathbb{R}^n,0)\to(\mathbb{R}^n,0)$$ making the following diagram commute

\begin{tikzpicture}[descr/.style={fill=white}]
\matrix(m)[matrix of math nodes,column sep={60pt,between origins},row
sep={40pt,between origins}]
{
&|[name=ka]| &|[name=kb]| (M,\sigma_f^{-1}(0)) &|[name=kc]| &|[name=kd]| (M',\sigma_g^{-1}(0)) &|[name=ke]|&|[name=kf]|   \\
&|[name=E]|   &|[name=F]| (\mathbb{R}^n,0) &|[name=G]|  &|[name=01]| (\mathbb{R}^n,0) &|[name=02]|\\
&|[name=A]|  &|[name=B]|  &|[name=C]| (\mathbb{R},0) &|[name=10]|  &|[name=20]| &|[name=30]|\\
};
\path[->,font=\scriptsize]
           (kb) edge node[auto]{\(\Psi\)}(kd)
          (F) edge node[auto]{\(\psi\)}(01)
          (kb) edge node[auto]{\(\sigma_f\)}(F)
          (kd) edge node[auto]{\(\sigma_g\)}(01)
          (F) edge node[auto]{\(f\)}(C)
          (01) edge node[auto]{\(g\)}(C)
          ;
\end{tikzpicture}
\end{enumerate}
We write $f\sim_{b-N} g$ if $f$ and $g$ are blow-Nash equivalent.
\end{definition}

That this is an equivalence relation follows from the following theorem.

\begin{theorem}[{\cite[Theorem 7.9]{Ca}}] Two germs of Nash functions 
$$f,g:(\mathbb{R}^n,0)\to(\mathbb{R},0)$$
are arc-analytically equivalent if and only if they are blow-Nash equivalent. 
\end{theorem}

\section{Invariance Under the Arc-Analytic Equivalence}\label{fourth part}
According to the result \cite[Theorem 4.9]{Fi} the virtual Poincaré polynomials of the motivic Milnor fibres are invariant under the arc-analytic equivalence. We shall begin by proving that this is also the case for the motivic Milnor tubes; the proof is in effect a minor modification of the proof of Fichou's result.

\begin{theorem}\label{Constancy}
If $f,g:(\mathbb{R}^n,0)\to(\mathbb{R},0)$ are arc-analytically equivalent germs of Nash functions then for any symbol $\epsilon\in\{\pm 1, <, >\}$ the virtual Poincaré polynomial of the respective motivic Milnor fibres and tubes $S^{\epsilon}(f)$ and $S^{\epsilon}(g)$, are equal. 
\end{theorem}
\begin{proof}
\begin{enumerate}
\item According to \cite[Proposition 7.9]{Ca} one has that $f\sim_{a} g$ if and only if $f\sim_{b-N} g$. Therefore there exist germs of Nash modifications 
$$\sigma_f: (M,\sigma_f^{-1}(0))\to(\mathbb{R}^n,0),\qquad \sigma_g: (M',\sigma_g^{-1}(0))\to(\mathbb{R}^n,0)$$
having the properties that the associated divisors of $f\circ\sigma_f$ and of $g\circ\sigma_g$ are normal crossings and that the associated divisors $\text{Jac}(\sigma_f)$ and $\text{Jac}(\sigma_g)$ are normal crossings, and that there exists a germ of Nash isomorphisms (\cite[Definition 3.2]{Fi})
$$\Psi: (M,\sigma_f^{-1}(0))\to(M',\sigma_g^{-1}(0))$$
preserving the multiplicities of the divisors $\text{div }\text{Jac}(\sigma_f)$ and $\text{div }\text{Jac}(\sigma_g)$ along the irreducible components of the corresponding exceptional divisors $\sigma_f^{-1}(0)$ and $\sigma_g^{-1}(0)$. 
\item\label{forbanna} Write 
$$(g\circ\sigma_g)^{-1}(0)=\bigcup_{i=1}^k E_i'$$ 
for the decomposition into irreducible components of $(g\circ\sigma_g)^{-1}(0)$ and let $p\in M'$ be a point. Since $\text{div}(g\circ\sigma_g)$ is a normal crossings divisor there exists an affine open set $U'\subset M'$ containing $p$ and a system of (étale) local coordinates $x_1,\dots, x_k$ on $U'\subset M'$ such that 
$$g\circ\sigma_g(x)=v(x)\prod_{i=1}^k x_i^{N_i}$$
for $v: U'\to \mathbb{R}$ an invertible Nash function. Since $\sigma_g^{-1}(0)$ is compact one can find a finite set $L\subset\mathbb{N}$ and open affine subsets $U_{l}', l\in L$ of a neighborhood of $\sigma_g^{-1}(0)$ such that 
$$\sigma_g^{-1}(0)\subset\bigcup_{l\in L} U_l'$$
For any subset $I\subset\{1,\dots,k\}$ consider the stratification of $(g\circ\sigma_g)^{-1}(0)$ with strata
$$E_I^{'o}=\bigcap_{i\in I} E_i'\setminus\bigcup_{j\notin I} E_j'.$$
If $E_I^{'0}\cap \sigma_g^{-1}(0)\neq \emptyset$ then it follows that 
$$E_I^{'o}=\bigcup_{l\in L} E_I^{'o}\cap U_l'.$$
\item Then $E_I^{'o}\in \mathcal{AS}$ since $g$ is a germ of Nash function, by \cite[Remark 6.10]{Ca}, and $U'\in\mathcal{AS}$ because $\mathcal{AS}$ is a constructible category which contains the category of Zariski constructibles, by (\cite[Theorem 4.4]{kurdyka}). Hence $E_I^{'o}\cap U'\in\mathcal{AS}$ again by (\cite[Theorem 4.4]{kurdyka}) because $\mathcal{AS}$ is a constructible category.\\

\item For each symbol $\epsilon\in\{\pm 1,<,>\}$ consider the classes in $K_0(\text{BSA})$ corresponding to $g$ and $U'$ as in Definition \hyperref[defn 12]{\ref*{defn 12}} and in the statement of \cite[Theorem 4.2]{CF}, namely
$$[R_{I, U'}^{\pm}(g)]=[(x,t)\in (E_I'^{o}\cap U')\times\mathbb{R}\ |\ t^{N_I} v(x)=\pm 1],$$
$$[R_{I, U'}^{>}(g)]=[(x,t)\in (E_I'^{o}\cap U')\times\mathbb{R}\ |\ t^{N_I} v(x)>0],$$
$$[R_{I, U'}^{<}(g)]=[(x,t)\in (E_I'^{o}\cap U')\times\mathbb{R}\ |\ t^{N_I} v(x)<0].$$
Since $\Psi$ is a Nash isomorphism the union of the preimages $E_I^{o}:=\Psi^{-1}(E_I'^{o})$ gives a stratification of 
$$(f\circ\sigma_f)^{-1}(0)$$
such that $\sigma_f^{-1}(0)$ is a union of strata. Moreover if $U'\subset M'$ is an affine open set then $U=\Psi^{-1}(U')$ is an affine open set because a Nash isomorphism is semialgebraic by \cite[Definition 3.2]{Fi}. Write 
$$\Psi=(\Psi_1,\dots,\Psi_n): (M, \sigma_f^{-1}(0))\to (M', \sigma_g^{-1}(0))$$ 
Since $\Psi$ is a germ of Nash isomorphism it preserves the multiplicities $N_i$ of $\text{div}(g\circ\sigma_g)$ along $E_i$ for each $i=1,\dots, k$ so commutativity of the diagram in the Definition \hyperref[defn 4.2]{\ref*{defn 4.2}} of the Blow-Nash equivalence yields 
$$f\circ\sigma_f(y)=g\circ\sigma_g\circ\Psi(y)= v(\Psi(y))\prod_{i=1}^k\Psi_i(y)^{N_i}$$
where $v\circ\Psi: U\to \mathbb{R}$ is an invertible Nash function. All said,
$$[R^{\pm}_{I, U}(f)]=[(y,t)\in (E_I^o\cap U)\times\mathbb{R}\ |\ t^{N_I} v(\Psi(y))=\pm 1],$$
$$[R^{>}_{I, U}(f)]=[(y,t)\in (E_I^o\cap U)\times\mathbb{R}\ |\ t^{N_I} v(\Psi(y))>0],$$
$$[R^{<}_{I, U}(f)]=[(y,t)\in (E_I^o\cap U)\times\mathbb{R}\ |\ t^{N_I} v(\Psi(y))<0]$$
are the classes in $K_0(\text{BSA})$ corresponding to $f$ and $U=\Psi^{-1}(U')$ as in the statement of \cite[Theorem 4.2]{CF}.

\item Consider now the realisations of these sets under the virtual Poincaré polynomial (Definition \hyperref[def polynome]{\ref*{def polynome}}). We claim that the underlying set $R^{\pm 1}_{I, U'}(g)$ is an intersection of $\mathcal{AS}$-sets. Indeed $E_I'^{o}\cap U'$ is $\mathcal{AS}$ by the above. Furthermore $v$ is Nash hence semialgebraic and real analytic so its graph is $\mathcal{AS}$, as in the second step of the proof of (\cite[Theorem 3.3]{Fi}). The inverse image of an $\mathcal{AS}$-set under an $\mathcal{AS}$-map is $\mathcal{AS}$, by (\cite[Theorem 4.4]{kurdyka}), hence $\{v(x)=\pm 1\}$ are $\mathcal{AS}$. Therefore $R^{\pm 1}_{I, U'}(g)$ is an intersection of $\mathcal{AS}$-sets hence is an $\mathcal{AS}$-set by (\cite[Theorem 4.4]{kurdyka}).

\item\label{skit} Again by the second step of the proof of \cite[Theorem 3.3]{Fi} the graph of $\Psi$ is $\mathcal{AS}$ so $\Psi$ is an $\mathcal{AS}$-map. Since it is an analytic isomorphism it is injective so the images of each of the sets $R_{I, U'}^{\pm 1}(g)$ under the injective $\mathcal{AS}$-map 
$$(x,t)\mapsto(\Psi^{-1}(x),t)$$
are $\mathcal{AS}$-sets, by (\cite[Theorem 4.4]{kurdyka}). Hence $R_{I, U}^{\pm 1}(f)$ are $\mathcal{AS}$-sets. Since $\Psi$ is a Nash isomorphism the application $(y,t)\mapsto(\Psi(y),t)$ gives a Nash isomorphism between $R^{\epsilon}_{I, U'}(g)$ and $R^{\epsilon}_{I, U}(f)$. Since $\beta_{\text{Var}}$ is invariant under Nash isomorphisms by \cite[Theorem 3.3]{Fi} it follows that
$$\beta(R_{I, U}^{\epsilon}(f))=\beta_{\text{Var}}(R_{I, U}^{\epsilon}(f))=\beta_{\text{Var}}(R_{I, \Psi(U)}^{\epsilon}(g))=\beta(R_{I, \Psi(U)}^{\epsilon}(g))$$
for any $I\subset\{1,\dots,k\}$ and any symbol $\epsilon\in\{\pm 1\}$. If $\bigcup_{l\in L} U_l'$ is a covering of a neighborhood of $\sigma_g^{-1}(0)$ by affine open sets as in the previous step \hyperref[forbanna]{\ref*{forbanna}} then 
$$\bigcup_{l\in L} \Psi^{-1}(U_l)$$
is a covering of a neighborhood of $\sigma_f^{-1}(0)$ by affine open sets hence 
$$\beta(E_{I}^{o,\epsilon}(f))=\sum_{S\subset L} (-1)^{|S|+1}\beta(R_{I, \cap_{s\in S} \Psi^{-1}(U_s)}(f))$$
$$=\sum_{S\subset L} (-1)^{|S|+1}\beta(R_{I, \cap_{s\in S} U_s}(g))=\beta(E_{I}^{o,\epsilon}(g))$$
so $\beta(S^{\epsilon}_f)=\beta(S^{\epsilon}_g)$ by \cite[Theorem 4.2]{CF}.
 
\item We shall prove that $\beta(S^{>}_f)=\beta(S^{>}_g)$, the proof of $\beta(S^{<}_f)=\beta(S^{<}_g)$ is similar. It follows by the Definition \hyperref[def polynome]{\ref*{def polynome}} of the virtual Poincaré polynomial on basic semialgebraic formulas that
$$\beta(R_{I, U}^{>}(f))=\frac{1}{4}\beta_{Var}([(y,t)\in E_I^o\cap U\times\mathbb{R}\ |\ t^{N_I}v (\Psi(y))=z^2])$$
$$-\frac{1}{4}\beta_{Var}([(y,t)\in E_I^o\cap U\times\mathbb{R}\ |\ t^{N_I}v(\Psi(y))=-z^2])$$
$$+\frac{1}{2}\beta_{Var}([(y,t)\in E_I^o\cap U\times\mathbb{R}\ |\ t^{N_I}v(\Psi(y))\neq 0]).$$
To show $\beta(R^{>}_{I, U}(f))=\beta(R^{>}_{I, \Psi(U)}(g))$ it therefore suffices to show that
$$\beta_{Var}([(x,t)\in E_I'^o\cap U'\times\mathbb{R}\ |\ t^{N_I}v(x)=\pm z^2])=$$
$$\beta_{Var}([(y,t)\in E_I^o\cap U\times\mathbb{R}\ |\ t^{N_I}v(\Psi(y))=\pm z^2])$$
and that
$$\beta_{Var}([(x,t)\in E_I'^o\cap U'\times\mathbb{R}\ |\ t^{N_I}v(x)\neq 0])=$$
$$\beta_{Var}([(y,t)\in E_I^o\cap U\times\mathbb{R}\ |\ t^{N_I}v(\Psi(y))\neq 0])$$
The underlying sets of these classes are in each of the above cases intersections of $\mathcal{AS}$-sets and Zariski constructible sets hence are $\mathcal{AS}$-sets by  (\cite[Theorem 4.4]{kurdyka}). Since $\Psi$ is a Nash isomorphism it follows that $(y,t,z)\mapsto(\Psi(y),t,z)=(x,t,z)$ gives a Nash isomorphism between
$$\{(y,t,z)\in E_I^o\cap U\times\mathbb{R}\times\mathbb{R}\ |\ t^{N_I} v(\Psi(y))=z^2\}$$
and
$$\{(x,t,z)\in E_I'^o\cap U'\times\mathbb{R}\times\mathbb{R}\ |\ t^{N_I} v(x)=z^2\}$$
and between 
$$\{(y,t,z)\in E_I^o\cap U\times\mathbb{R}\times\mathbb{R}\ |\ t^{N_I}v(\Psi(y))\neq 0\}$$
and 
$$\{(x,t,z)\in E_I'^o\cap U'\times\mathbb{R}\times\mathbb{R}\ |\ t^{N_I} v(x)\neq 0\}.$$
Thus $\beta(R^{>}_{I, U}(f))=\beta(R^{>}_{I, \Psi(U)}(g))$ for all $I\subset \{1,\dots, k\}$ by \cite[Theorem 3.3]{Fi} so the same argument as in the previous step \hyperref[skit]{\ref*{skit}} gives
$$\beta(S^{>}(f))=\beta(S^{>}(g)).$$
\end{enumerate}
\end{proof}

\section{Arc-Wise Analytic Triviality}\label{fifth part}
\subsection{Arc-Wise Analytic Triviality of Function Germs}
Before continuing we need to recall the notion of arc-wise analytic triviality for germs of families of analytic functions; we refer to the article \cite{PP} for further details. In the following we shall denote by $\mathbb{K}$ either the field of complex, or real numbers.
\begin{definition}[{\cite[§ 8.3]{PP}}] Let $T\subset\mathbb{K}^m$ be an open analytic set. A $\mathbb{K}$-analytic family $f: T\times(\mathbb{K}^{n-1},0)\to\mathbb{K}$ of germs of $\mathbb{K}$-analytic functions is \emph{arc-wise analytically trivial} over $T$ if there exists a neighborhood $\Lambda\subset\mathbb{K}^m\times\mathbb{K}^{n-1}$ of $T\times\{0\}$ and a neighborhood $\Lambda_0\subset\mathbb{K}^{n-1}$ of $\{0\}$ and an arc-wise analytic trivialisation \cite[Definition 1.2]{PP}
$$\sigma: T\times\Lambda_0\to\Lambda$$
such that $f\circ \sigma (t,x)=f(0,x)$.
\end{definition}
The key result which is going to be used later on is the following theorem. Here we only assume that $T\subset\mathbb{K}^m$ is analytic.

\begin{theorem}[{\cite[Theorem 8.5]{PP}}]\label{theorem 8.5} Let $f: T\times(\mathbb{K}^{n-1},0)\to\mathbb{K}$ be a $\mathbb{K}$-analytic family of $\mathbb{K}$-analytic function germs and let $t_0\in T$. Then there exists a neighborhood $U\subset T$ of $t_0$ and a $\mathbb{K}$-analytic subset $Z\subset U$ with $\dim Z<\dim U$ such that $f$ is arc-wise analytically trivial along $U\setminus Z$.
\end{theorem}
\subsection{The Main Theorem}
The main result of this chapter is the following Theorem \hyperref[main theorem family]{\ref*{main theorem family}}. The content of the proof is to show that in the proof of the Theorem \hyperref[theorem 8.5]{\ref*{theorem 8.5}} above, the constructed arc-analytic trivialisation of $\sigma$ has the property that the jacobian determinant of $\sigma_t:=\sigma(t, \cdot)$ is bounded where it is defined and thus provides the wanted arc-analytic equivalence. 

\begin{theorem}\label{main theorem family} Let $T\subset\mathbb{R}^m$ be a real algebraic set and let $f: T\times(\mathbb{R}^n,0)\to\mathbb{R}$ be a polynomial function of polynomial map germs such that $f(t,0)=0$ for any $t\in T$. There exists a locally finite real analytic stratification $\mathcal{S}$ of $T$ such that if $S\in\mathcal{S}$ is a stratum then
$$f_t=f(t, \cdot)\sim_a f_{t'}=f(t,\cdot),\qquad\forall t,t'\in S$$
in the notation of Definition \hyperref[arc equivalence]{\ref*{arc equivalence}}. In particular 
$$\beta(S^{\epsilon}(f_t))=\beta(S^{\epsilon}(f_{t'})),\qquad \forall t,t'\in S\in\mathcal{S},\qquad \forall \epsilon\in\{\pm 1, <, >\}.$$ 
Moreover if $T$ is compact then $\mathcal{S}$ is a finite stratification.
\end{theorem}
\begin{proof}
\begin{enumerate}
\item Recalling the proof of the Theorem \hyperref[theorem 8.5]{\ref*{theorem 8.5}} (\cite[Theorem 8.5]{PP}) let us write 
$$F: T\times\mathbb{R}\times\mathbb{R}^n\to\mathbb{R},\qquad F(t,z, x)=z-f(t,x)$$
where $x=(x_1,\dots, x_n)\in \mathbb{R}^n$. Define    
$$Y:=\{F=0\}\subset T\times\mathbb{R}\times\mathbb{R}^n.$$
Fix $t_0\in T$. Then according to \cite[Lemma 6.3]{PP} there exists a neighborhood $U(t_0)$ containing $t_0$, a proper analytic subset 
$$Z(t_0)\subset U(t_0),\qquad \dim Z(t_0)<\dim U(t_0)$$
and a transverse local system $(F_i)_{i=0}^{n+1}$ of Zariski equisingular pseudo-polynomials (see \cite[Definition 4.1, 3.1]{PP}) such that $F_{n+1}$ is the Weierstrass polynomial of $F$ in $(t,0)$ for all $t\in U(t_0)\setminus Z(t_0)$. By \cite[Theorem 3.3]{PP} there exist for any such local system a real number $\epsilon>0$ and neighborhoods 
$$\mathbb{B}_{\epsilon}=\{t\in U(t_0)\setminus Z(t_0)\ |\ \lVert t-t_0\rVert\leq\epsilon\}\subset\mathbb{R}^m,$$
$$\Omega_0\subset\mathbb{R}^{n+1},\qquad \Omega\subset\mathbb{R}^{n+m+1}$$
of $t_0$ and of the respective origins and a homeomorphism 
$$\Phi: \mathbb{B}_{\epsilon}\times\Omega_0\to\Omega$$
preserving $\{z_0=f(t,x)\}$ for each $z_0\in\mathbb{R}$ (see \cite[Definition 1.2]{PP}) such that $(z_0, x)\in \Omega_0$, and such that the conditions $(Z2-Z3)$ of \cite[Theorem 3.3]{PP} are satisfied, that is for all $t\in \mathbb{B}_{\epsilon}$ and all $(z, x_1,\dots, x_n)\in \Omega_0$ the following holds,\\
\begin{enumerate}
\item\label{Z zero} $\Phi$ is an arc-analytic trivialisation of the standard projection $\Omega\to \mathbb{B}_{\epsilon}$. In particular $\Phi$ is real analytic in $t$ and arc-analytic in $(z,x_1,\dots, x_n)$ and $\Phi^{-1}$ is arc-analytic.
\item\label{Z one} $\Phi(t, 0)=(t, 0)$ and 
$$\Phi(t_0, z, x_1,\dots, x_n)=(t_0, z, x_1,\dots, x_n).$$
\item\label{Z two} $\Phi(t,z, x_1,\dots,x_n)=(t, \Psi_1(t,z),\dots,\Psi_{n+1}(t,z,x_1,\dots,x_n))$.
\item\label{Z three} For each $i=2,\dots,n+1$,  
$$\Psi_i(t, z, x_1,\dots, x_{i-1},\bullet):\mathbb{R}\to\mathbb{R}$$
are bi-Lipschitz with the Lipschitz constants of $\Psi_i$ and of $\Psi_i^{-1}$ independent on $(t, z, x_1,\dots, x_{i-1})$ and
$$\Psi_1(t, \cdot): \mathbb{R}\to \mathbb{R}$$
is bi-Lipschitz with Lipschitz constants of $\Psi_1$ and of $\Psi_1^{-1}$ independent of $t$.\\
\end{enumerate}
Write $\Phi(t,z,x)=(t,\Psi(t,z,x))$ (so that $\Psi: \mathbb{B}_{\epsilon}\times \Omega_0\to \Omega_0$) and consider the functions
$$\Psi_t=\Psi(t, \cdot,\cdot): \Omega_0\to \Omega_0,\qquad t\in U(t_0)\setminus Z(t_0)$$
$$\Phi_t=(t, \Psi_t): \Omega_0\to \mathbb{B}_{\epsilon}\times \Omega_0$$
The above conditions (\hyperref[Z two]{\ref*{Z two}}) and (\hyperref[Z three]{\ref*{Z three}}) imply, as in the proof of \cite[Proposition 3.6]{PP}, that there exist constants $\tilde{C},\tilde{c}>0$, which are by condition (\hyperref[Z three]{\ref*{Z three}}) independent of $t$, such that the Jacobian determinants are bounded
\begin{equation}\label{import}
\tilde{c}\leq|\det\text{Jac}\Psi_t|\leq \tilde{C}
\end{equation}
there where they are defined. In particular $\Psi_t$ is a diffeomorphism for $t\in U(t_0)\setminus Z(t_0)$ (see also \cite[Proposition 3.5]{PP}). Denote by  
$$\pi: \mathbb{R}^m\times\mathbb{R}\times\mathbb{R}^n\to \mathbb{R}^m\times\mathbb{R}^n,\qquad \pi(t,z,x)=(t,x)$$ 
the standard projection and put
$$\tilde{\pi}: \mathbb{R}^{n+1}\to\mathbb{R}^n,\qquad \tilde{\pi}(z,x)=x$$
Let $\Lambda_0=\tilde{\pi}(\Omega_0)$ and $\Lambda=\pi(\Omega)$. Then 
$$\sigma: U(t_0)\setminus Z(t_0)\times \Lambda_0\to\Lambda$$
given by
$$\sigma(t,x)=\pi\circ\Phi(t, f_0(x),x)$$
is an arc-wise analytic trivialisation of $f$ along $U(t_0)\setminus Z(t_0)$, by \hyperref[theorem 8.5]{\ref*{theorem 8.5}} (\cite[Theorem 8.5]{PP}). 
\item Consider
$$\sigma_t=\sigma(t,\cdot): \Lambda_0\to \Lambda.$$
Then by construction,
$$\sigma_t=\pi\circ \Phi(t, f_0(\cdot), \cdot)=$$
$$=\pi\circ (t, \Psi(t, f_0(\cdot), \cdot)=\tilde{\pi}\circ \Psi_t(f_0(\cdot),\cdot)$$
so that one can write
$$\sigma_t=\tilde{\pi}\circ \Psi_{t| H}=\tilde{\pi}\circ \Psi_t\circ g$$
where $H\subset \Omega_0$ is the smooth hypersurface $H=\Gamma(f_0)\cap \Omega_0$ where 
$$\Gamma(f_0)=\{(z, x)\in \mathbb{R}^{n+1}\ |\ z=f_0(x)\}$$
and where 
$$g: \Lambda_0\to \Omega_0,\qquad g(x)=(f_0(x),x).$$
\item\label{rovhar} Let $t\in U(t_0)\setminus Z(t_0)$. Then 
$$\text{Jac}(\sigma_t)(x)=\text{Jac}(\tilde{\pi})(\Psi_t(f_0(x),x))\text{Jac}(\Psi_t)(f_0(x),x)\text{Jac}(g)(x)$$
by the chain rule. Taking local coordinates one finds that the Jacobian matrix of $g$ has the form 
\[
\text{Jac}(g)(x)= \begin{pmatrix}
(\partial f_0/ \partial x_i)(x) & I_n
\end{pmatrix}
\]
with $I_n=\text{Id}_{\mathbb{R}^n}$ the identity matrix, so  $\text{Jac}(g)$ has full rank. By \hyperref[import]{\ref*{import}} one has that $\text{Jac}(\Psi_{t})$ likewise has full rank, since its determinant is nonzero. 

\item\label{morrhar} We claim that, after replacing if necessary $U(t_0)$ by a smaller neighborhood of the origin, for all $t\in U(t_0)\setminus Z(t_0)$ the restriction $\tilde{\pi}_{| \Psi_t(H)}$ of $\tilde{\pi}$ to $\Psi_t(H)$ is a submersion. To show this, remark that
$$\tilde{\pi}_{|H}: H\to \mathbb{R}^n$$
is a submersion, by the previous step (\hyperref[rovhar]{\ref*{rovhar}}). This implies that there exists a real number $\delta>0$ such that if $\mathbb{B}_{\delta}\subset \mathbb{R}^{n+1}$ denotes the open ball of radius $\delta$ centered at the origin then $\bar{\mathbb{B}}_{\delta}\subset \Omega_0$ and 
$$\forall s\in \tilde{\pi}(H\cap \partial \bar{\mathbb{B}}_{\delta}),\qquad \tilde{\pi}^{-1}(s)\pitchfork (H\cap \partial\bar{\mathbb{B}}_{\delta}).$$
In particular
$$\tilde{\pi}_{|H\cap \partial \bar{\mathbb{B}}_{\delta}}: H\cap \partial \bar{\mathbb{B}}_{\delta}\to \mathbb{R}^n$$
is a submersion. Replacing $\Omega_0$ by $\mathbb{B}_{\delta}$ we can therefore assume that 
$$\tilde{\pi}_{| H\cap \bar{\Omega}_0}: H\cap \bar{\Omega}_0\to \mathbb{R}^n$$
is a submersion. Remark that $t\mapsto \Psi_t$ is continuous by (\hyperref[Z zero]{\ref*{Z zero}}) and that $\Psi_{t_0}$ is the identity map on $\Omega_0$ by (\hyperref[Z one]{\ref*{Z one}}). Since being a submersion is an open condition and since $H\cap \bar{\Omega}_0$ is compact it follows that there exists an open neighborhood $U'(t_0)\subset \mathbb{R}^m$ of $t_0$ such that $U'(t_0)\subset U(t_0)$ and such that 
$$\tilde{\pi}_{|\Psi_t(H)}: \Psi_t(H)\to \mathbb{R}^n$$
is a submersion. This proves the claim.

\item Fix $t\in U(t_0)\setminus Z(t_0)$. By the previous step (\hyperref[morrhar]{\ref*{morrhar}}) $\det\text{Jac}(\sigma_t)$ is nonzero at any $x\in \Lambda_0$. Recall Hadamard's determinant inequality: if $A\in M(k,k, R)$ is a $k\times k$ matrix over a commutative ring $R$ then
$$|\det(A)|\leq \lVert A\cdot e_1\rVert\dots \lVert A\cdot e_k\rVert$$
where $A\cdot e_i$ are the column vectors of $A$. Again by (\hyperref[Z three]{\ref*{Z three}}),
$$\lVert \text{Jac}\Psi_t\cdot e_i\rVert\leq K,\qquad i=1,\dots, n+1.$$
for a constant $K$ not depending on $t$. Let 
$$M:=\text{max}_{1\leq i\leq n}\text{sup}_{x\in\bar{\Lambda}_0} \partial f_0/\partial x_i (x)$$
where $\bar{\Lambda}_0\subset \mathbb{R}^n$ is the closure of $\Lambda_0$ for the euclidean topology. Since $\text{Jac}(\sigma_t)$ is invertible, the chain rule and Hadamard's inequality yield
$$|\det\text{Jac}(\sigma_t^{-1})(\sigma_t(x))|=\frac{1}{|\det\text{Jac}(\tilde{\pi}\circ \Psi_{t|H})(x)|}$$
$$\geq \frac{1}{\prod_{i=1}^{n} \lVert \text{Jac}(\tilde{\pi}\circ \Psi_{t|H})\cdot e_i\rVert}$$
$$\geq \frac{1}{l M \prod_{i=1}^{n+1} \lVert \text{Jac}(\Psi_t) \cdot e_i\rVert } \geq \frac{1}{l M K^{n+1}}$$
for some nonzero $l\in\mathbb{N}$. Therefore there exists a positive constant\\
$c=1/(l M K^{n+1})$ not depending on the parameter $t$ (because $K$ and $M$ are independent of $t$), such that 
$$|\det\text{Jac}(\sigma_t^{-1})|\geq c.$$ 
By the definition \cite[Definition 1.2]{PP} of an arc-wise trivialisation $\Phi_t$ is arc-analytic and its inverse is arc-analytic hence $\sigma_t$ is arc-analytic with arc-analytic inverse. And as its determinant is bounded where it is defined it follows by Definition \hyperref[arc equivalence]{\ref*{arc equivalence}} of the arc-analytic equivalence that for each $t\in U(t_0)\setminus Z(t_0)$ the function $\sigma_t^{-1}$ yields an arc-analytic equivalence $f_t\sim_a f_0$. Therefore, 
$$\beta(S^{\epsilon}(f_t))=\beta(S^{\epsilon}(f_{t_0})),\qquad\forall t\in U(t_0)\setminus Z(t_0).$$
by Theorem \hyperref[Constancy]{\ref*{Constancy}}.
\item For each $t_0\in T$ one finds a finite analytic stratification of a neighborhood of $t_0$ as follows. Let $U_0=U(t_0)$ be the neighborhood of $t_0$ constructed in the previous step and let $Z_0=Z(t_0)$ be the corresponding proper analytic subset. Repeating the previous step we find for each $t_1\in Z_0$ a neighborhood $U_1\subset Z_0$ of $t_1$, a proper analytic subset $Z_1\subset U_1$ with $\dim Z_1<\dim U_1\leq \dim Z_0$ and an arc-wise analytic trivialisation
$$\sigma_1: U_1\setminus Z_1\times \Lambda_0^1\to\Lambda^1,\qquad f\circ\sigma_1(t,x)=f(0,x).$$
By induction on the dimension one finds a finite list of subsets 
$$S_0(t_0),\dots, S_n(t_0)\subset U_0$$
giving a stratification $\mathcal{S}(t_0)$ of the neighborhood $U_0$ of $t_0\in T$ such that 
$$\beta(S^{\epsilon}(f_t))=\beta(S^{\epsilon}(f_{t'})),\qquad\forall t,t'\in S_i(t_0)$$

\item One takes a covering of $T$ by neighborhoods of the form $U=U(t)$ and one considers for each $t\in T$ the corresponding finite stratifications $\mathcal{S}(t)$. One defines $\mathcal{S}$ as a common refinement of the stratifications $\mathcal{S}(t)$. If $T$ is compact then there exists a finite subcovering 
$$T\subset\bigcup_{i=1}^l U(t_i)$$
with finite stratifications $\mathcal{S}(t_i)$ of $U(t_i)$, for each $i=1,\dots,l$ so in this case the common refinement $\mathcal{S}$ is finite as well. 
\end{enumerate}

\end{proof}

\bibliographystyle{plain}
\bibliography{biblio.bib}

\end{document}